\documentclass[a4paper,12pt,leqno]{article}
\usepackage{latexsym}
\usepackage[all]{xy}

\usepackage{amssymb} 
\usepackage{amsmath} 
\usepackage{theorem}

\def\Z{{\mathbb{Z}}}
\def\K{{\mathbb{K}}}
\def\R{{\mathbb{R}}}
\def\A{{\mathcal{A}}}
\def\B{{\mathcal{B}}}

\DeclareMathOperator{\codim}{codim}

\DeclareMathOperator{\Der}{Der}

\DeclareMathOperator{\Shi}{Shi}

\numberwithin{equation}{section}

\newcommand{\owari}{\hfill$\square$}

\theoremstyle{break}
\newtheorem{theorem}{Theorem}[section]
\newtheorem{prop}[theorem]{Proposition}
\newtheorem{cor}[theorem]{Corollary}
\newtheorem{lemma}[theorem]{Lemma}
\newtheorem{define}[theorem]{Definition}
\newtheorem{definetheorem}[theorem]{Theorem-Definition}
\newtheorem{rem}[theorem]{Remark}

\newtheorem{example}[theorem]{Example}
\newtheorem{problem}[theorem]{Problem}
\newtheorem{conj}[theorem]{Conjecture}

\title{Divisionally free arrangements of hyperplanes
}
\author{
Takuro Abe
\footnote
{
Institute of Mathematics for Industry, 
Kyushu University,
Fukuoka 819-0395, Japan.
Email:abe@imi.kyushu-u.ac.jp\ 
\textit{2010 Mathematics Subject Classification}. 32S22, 52S35.
}
}
\date{\today}

\begin{document}

\maketitle

\begin{abstract}
Let $(\A,\A',\A^H)$ be the triple of hyperplane arrangements. 
We show that the freeness of $\A^H$ 
and the division of $\chi(\A;t)$ by $\chi(\A^H;t)$ imply the freeness of $\A$. 
This ``division theorem'' 
improves the famous addition-deletion theorem, and it has  
several applications, which include 
a definition of ``divisionally 
free arrangements''. It is a strictly larger class of free arrangements 
than the classical important class of 
inductively free arrangements. 
Also, whether an arrangement is divisionally free or not is 
determined by the combinatorics. 
Moreover, we show that 
a lot of recursively free arrangements, to 
which almost all known free arrangements 
are belonging, are divisionally free. 
\end{abstract}

\section{Main results}
Let $V$ be an $\ell$-dimensional vector space over an arbitrary field $\K$ with $\ell \ge 1$, 
$S=\mbox{Sym}(V^*)=\K[x_1,\ldots,x_\ell]$ its coordinate ring and 
$\Der S:=\oplus_{i=1}^\ell S \partial_{x_i}$ the module of $\K$-linear $S$-derivations. 
A \textbf{hyperplane arrangement} $\A$ is a finite set of hyperplanes in $V$. We say that $\A$ is 
\textbf{central} if every hyperplane is linear. In this article every arrangement is 
central unless otherwise specified. In the central cases, we fix a linear form $\alpha_H \in V^*$ 
such that $\mbox{ker}(\alpha_H)=H$ for each $H \in \A$. 
An \textbf{$\ell$-arrangement} is an arrangement in an $\ell$-dimensional vector space.
Let $L(\A):=
\{\cap_{H \in \B} H \mid \B \subset \A\}$ be an \textbf{intersection lattice}. $L(\A)$ has a 
partial order by reverse inclusion, which equips $L(\A)$ with a poset structure. 
For $X \in L(\A)$, define the \textbf{localization} $\A_X$ of $\A$ at $X$ by 
$\A_X:=\{H \in \A \mid H \supset X\}$, which is a subarrangement of $\A$. Also, the 
\textbf{restriction} $\A^X$ of $\A$ onto $X$ is defined by 
$\A^X:=\{H \cap X \mid H \in \A \setminus \A_X\}$, which is an arrangement in 
$X \simeq \K^{\dim X}$. Let $L_i(\A):=
\{X \in L(\A) \mid \codim_V X=i\}$. Also, 
we use some notations in \S2.

In the study of hyperplane arrangements, its algebraic structure $D(\A)$ is well-studied. 
The \textbf{logarithmic derivation module} $D(\A)$ is defined by 
$$
D(\A):=\{\theta \in \Der S \mid \theta(\alpha_H) \in S \cdot \alpha_H \ (\forall H \in \A)\}.
$$
We say that $\A$ is \textbf{free} with \textbf{exponents} $\exp(\A)=(d_1,\ldots, d_\ell)$ if 
$D(\A)$ is generated as an $S$-module by $S$-independent homogeneous generators 
$\theta_1,\ldots,\theta_\ell$ with $\deg \theta_i=d_i\ (i=1,\ldots,\ell)$. The study of free arrangements was initiated by Terao in \cite{T1}, 
and has been playing the central role in this area.
Recently, there have been 
several studies to determine when 
$\A$ is free, e.g., \cite{A}, \cite{AY}, \cite{Y1}, \cite{Y2} and so on. However, it is still very difficult to 
determine the freeness. Freeness of 
arrangements implies several interesting geometric and combinatorial properties of $\A$. For 
example, see \cite{A1}, \cite{A} and \cite{T2}. In particular, the most important 
result among them 
is Terao's factorization theorem (Theorem \ref{fact}) in \cite{T2}, which asserts that if $\A$ is free with 
$\exp(\A)=(d_1,\ldots,d_\ell)$, then its \textbf{characteristic polynomial} $\chi(\A;t)$ (essentially 
this is the same as the topological Poincar\`{e} polynomial $\pi(\A;t)$ of the 
complement $M(\A):=V \setminus \cup_{H \in \A} H$ of $\A$ in $V$ when $\K=\mathbb{C}$) 
factorizes into $\chi(\A;t)=\prod_{i=1}^\ell (t-d_i)$. 
When $\A \neq \emptyset$, it is known that 
$(t-1)$ divides $\chi(\A;t)$. Let $\chi_0(\A;t):=
\chi(\A;t)/(t-1)$ when $\A \neq \emptyset$.

The most useful method to construct free arrangements is the addition-deletion theorems 
(Theorem \ref{adddel})  
by Terao in \cite{T1}. Let us recall it. For a central arrangement $\A$ and $H \in \A$, define 
$
\A':=\A \setminus \{H\}.
$
Combining this with the restriction $\A^H$ of $\A$ onto $H \in \A$, 
we call $(\A,\A',\A^H)$ the \textbf{triple}. 
The addition-deletion theorem enables us to 
determine the freeness of all these three when we know the freeness of any two of 
them with some information on exponents. For example, the addition theorem 
asserts that if both $\A'$ and $\A^H$ are free with $\exp(\A') \supset \exp (\A^H)$, then 
$\A$ is also free with certain exponents. By the factorization theorem above, in this case, 
it follows that $\chi(\A^H;t)$ divides $\chi(\A';t)$ (the division of polynomials is often denoted by
$\chi(\A^H;t) \mid \chi(\A';t)$).  Since there is a famous  
deletion-restriction formula (see Theorem \ref{delres})
$$
\chi(\A;t)=\chi(\A';t)-\chi(\A^H;t),
$$
it holds that $\chi(\A^H;t)$ also divides $\chi(\A;t)$ in this case. So the addition-deletion 
theorem contains 
a statement about divisions of polynomials, but these divisions have not been studied so much. 

The aim of this article is to give a consideration on this aspect, i.e., the division of 
characteristic polynomials of these triples. The main result in this article is as follows.

\begin{theorem}[Division theorem]
Let $\A$ be a central $\ell$-arrangement. 
Assume that there is a 
hyperplane $H \in \A$ such that $\chi(\A^H;t)$ divides $\chi(\A;t)$ and that 
$\A^H$ is free. Then 
$\A$ is free.
\label{newadd}
\end{theorem}

Theorem \ref{newadd} has several advantages compared to previous results in 
both theoretical and practical senses. First, 
let us show a practical advantage of Theorem \ref{newadd} in the 
following example.

\begin{example}



Let us consider the arrangement $\B$ in $\R^4$ defined by 
$$
\B=
(\prod_{i=1}^4 x_i) \prod_{a_2,a_3,a_4 \in \{\pm 1\}} (x_1 +a_2 x_2+a_3 x_3 +a_4 x_4)=0.
$$
$\B$ is a restriction of the famous counterexample to Orlik's famous conjecture by 
Edelman and Reiner in \cite{ER} onto a coordinate hyperplane. $\B$ 
is known to be free. Here we can check it soon by 
Theorem \ref{newadd} combinatorially.
Let 
$$
\mathcal{C}:=\B^{x_4=0}=
(\prod_{i=1}^3 x_i) \prod_{a_2,a_3\in \{\pm 1\}} (x_1 +a_2 x_2+a_3 x_3)=0.
$$
Then it is easy to check that 
\begin{eqnarray}
\label{A}
\chi(\B;t)&=&(t-1)(t-3)^2(t-5),\\
\chi(\mathcal{C};t)&=&(t-1)(t-3)^2,\nonumber \\
\chi(\mathcal{C}^{x_3=0};t)&=&(t-1)(t-3).\nonumber
\end{eqnarray}
Since every $2$-central arrangment is free (e.g., see Lemma \ref{basic}),  
applying Theorem \ref{newadd} twice shows that $\B$ is free. 
Note that the freeness of $\B$ is checked by using only combinatorial information 
$L(\B)$.
\label{ERexample0}
\end{example}

There are several ways to understand Theorem \ref{newadd}. 
If we emphasize the aspect of a generalization of the addition-deletion theorem \ref{adddel}, 
Theorem \ref{newadd} can be 
formulated as Theorems \ref{adddiv} and 
\ref{same}. Also, it can be regarded as a modified converse of 
Orlik's conjecture, see Remark \ref{Orlik}. Here we focus on the relation between 
the freeness and combinatorics.

As in Example \ref{ERexample0}, 
Theorem \ref{newadd} allows us to 
check the freeness of $\A$ by constructing a divisional ``tower'' of 
characteristic polynomials as in (\ref{A}). In other words, we obtain the following 
completely combinatorial algorithm to show freeness:


\begin{theorem}
An $\ell$-arrangement $\A$ is free if 
there is a flag 
$$
V=X_0 \supset X_1 \supset \cdots \supset X_{\ell-2}
$$
such that 
$X_i \in L_i(\A)$ for 
$i=0,\ldots,\ell-2$ and 
$\chi(\A^{X_{i+1}};t) \mid \chi(\A^{X_{i}};t)$ for 
$i=0,\ldots,\ell-3$.
Let us call this flag a \textbf{divisional flag}.
\label{tower}
\end{theorem}



Apparently, whether $\A$ has a divisional flag or not is purely combinatorial.
Hence Theorem \ref{tower} 
gives us purely combinatorial method to check the freeness.
One of these theoretical applications is 
%
to investigate the relation between freeness and combinatorics of 
the famous Shi arrangements
(see \S 6 for 
the notation):

\begin{theorem}
The extended Shi arrangement has a divisional flag. Hence 
the freeness of Shi arrangement depends only on its intersection lattice.
\label{Shi}
\end{theorem}

To determine the relation between freeness and combinatorics has been 
an imporant but very hard problem. However, in this article, we can show the similar results as in Theorem \ref{Shi} due to Theorem \ref{newadd} (e.g., 
Example \ref{31} and Theorem \ref{akr}), which is very useful for that purpose.
To make this framework systematically, we may introduce a new class 
of free arrangements, called the class of 
\textbf{divisionally free arrangements}.

\begin{define}[Divisionally free arrangements]
We say that an $\ell$-arrangement $\A$ is 
\textbf{divisionally free} if $\A$ has a divisional flag
in 
Theorem \ref{tower}, which is combinatorially determined by $L(\A)$. 
The set $\mathcal{DF}_\ell$ consists of 
all the divisionally free $\ell$-arrangements, and $\mathcal{DF}:=
\cup_{\ell \ge 1} \mathcal{DF}_\ell$.
\label{DF0}
\end{define}

For details on divisionally free arrangements, see Theorem \ref{DFTC}. The most 
famous and important class of free arrangements 
with the same properties is the class of inductively free arrangements introduced by 
Terao in \cite{T1} (see Definition \ref{IF}). 
Let $\mathcal{IF}$ denote the class of inductively free 
arrangements. 
In fact, the class of divisionally free arrangements is strictly larger than 
that of inductively free arrangements.

\begin{theorem}
$
\mathcal{IF} \subsetneq \mathcal{DF}$.
\label{IFDF}
\end{theorem}

To prove Theorem \ref{IFDF}, recent developments on 
the freeness of unitary reflection arrangements (e.g., 
\cite{AHR}, \cite{HR}) play the key roles.
Hence Theorem \ref{newadd} develops 
the theory of free arrangements from the viewpoint of 
Terao's conjecture (Conjecture \ref{TC}) by strictly 
enlarging the inductively free arrangements. Also, we can show that 
a lot of recursively free 
arrangements are divisionally free in Theorem \ref{nb}. Moreover, 
as in Theorem \ref{Shi}, applications to those related to root systems 
are also given. 
Hence Theorem \ref{newadd} and 
the divisionally free arrangements 
give rise to an essential progress to the theory of free 
arrangements and its combinatorics.

To prove Theorem \ref{newadd}, the key ingredient is the following fact that, the division of 
characteristic polynomial commutes with localizations along 
the hyperplane in codimension-three :

\begin{theorem}[Localization and Remainder Theorem]
Let $\A$ be a central $\ell$-arrangement and $H \in \A$. Let us consider the 
polynomial division 
$$
\chi_0(\A;t)=(t-(|\A|-|\A^H|))\chi_0(\A^H;t)+r(t) 
$$
with the remainder $r(t)=\sum_{i=0}^{\ell-3} (-1)^i  r_i t^{\ell-3-i}$. 
Then $r_0 \ge 0$. Moreover, if 
$r_0=0$, then $\chi(\A_X^H;t)$ divides $\chi(\A_X;t)$ for all $X \in L_2(\A^H)$, and 
$\A$ is locally free along $H$ in codimension three.
In particular, if $\chi(\A^H;t)$ divides $\chi(\A;t)$ for a hyperplane $H \in \A$, then the 
same statement holds.
\label{codim3}
\end{theorem}

The statement in Theorem \ref{codim3} is non-trivial. For example, even when 
$\A$ is free, that statement does not hold in general, e.g, see Example \ref{nondiv}. 
Also, a certain converse of Theorem \ref{codim3} holds, see Theorem \ref{tower2}. 

If $\chi(\A^H;t)$ divides $\chi(\A;t)$, then clearly 
$r_0=0$. Hence the first part of Theorem \ref{codim3} implies the second one.
The 
meaning of $r_0 \ge 0$ in Theorem \ref{codim3} can be seen in Remark \ref{r0}.

\begin{rem}
As in Theorem \ref{codim3}, 
the condition $\chi(\A^H;t) \mid \chi(\A;t)$ in Theorem \ref{newadd} can be replaced by 
$r_0=0$ in terms of Theorem \ref{codim3}. 
\label{r0division}
\end{rem}

Following Remark \ref{r0division}, we can show the following 
equivalent condition to divisionally free arrangements. 

\begin{theorem}
The following three conditions are equivalent:
\begin{itemize}
\item[(1)]
$\A$ is divisionally free.
\item[(2)]
There is a flag $V=X_0 \supset X_1 \supset \cdots \supset X_{\ell-2}$ such 
that $X_i \in L_i(\A)\ (i=0,\ldots,\ell-2)$ and 
that $\chi_0(\A;t)=(t-(|\A^{X_{\ell-2}}|-1))\prod_{i=0}^{\ell-3} (t-(|\A^{X_i}|-|\A^{X_{i+1}}|))$.
\item[(3)]
There is a flag $V=X_0 \supset X_1 \supset \cdots \supset X_{\ell-2}$ such 
that $X_i \in L_i(\A)\ (i=0,\ldots,\ell-2)$ and 
that 
$$
b_2(d\A)=\sum_{i=0}^{\ell-3} 
(|\A^{X_i}|-|\A^{X_{i+1}}|)(|\A^{X_{i+1}}|-1),
$$
where $b_2(d\A)$ is the coefficient of $t^{\ell-3}$ in $\chi_0(\A;t)$.
\end{itemize}
\label{b2free1}
\end{theorem}

Hence to determine whether $\A \in \mathcal{DF}$, it suffices to compute 
$b_2(d\A)$, the second Betti number of the complement of $\A$ when $\K=\mathbb{C}$, and 
the number of restricted hyperplanes $|\A^{X_i}|$, which are numerical and combinatorial. 
We hope that Theorem \ref{b2free1} (3) enables us to 
make an efficient computer program to check the divisional freeness.

Though Theorem \ref{newadd} is a generalization of the addition-deletion theorems, 
the proofs of Theorems \ref{newadd} and \ref{codim3} deeply depend on 
the new development of the theory of multiarrangements, i.e., we use
the defintion of multiarrangements by Ziegler in \cite{Z}, 
Yoshinaga's criterion in \cite{Y1} and its refinement in \cite{AY}, 
development of basic tools to treat them in \cite{ATW} and \cite{ATW2}, 
and the same statement as Theorem \ref{newadd} when $\ell=3$ in \cite{A}. It is interesting 
to see that  
in the statement of Theorem \ref{newadd} there are no multiarrangements.

The organization of this article is as follows. 
In \S2 we introduce several results used for the proof of results in \S1. 
In \S3 we prove our main theorems in \S1. 
In \S 4 we investigate  
the most important application; the definition of divisionally free arrangements in 
Definitions \ref{DF0} and \ref{DF}. 
We show that divisionally free arrangements contain 
all the inductively free arrangements, and
whether an arrangement is divisionally free or not is determined by 
only combinatorics.
Moreover, 
there is a divisionally free arrangement which is not inductively free.
In \S 5 we consider the relation between divisionally and recursively free 
arrangements to which almost all known free arrangments belong. 
In \S 6 we give applications to the arrangements related to root systems. In particular, 
we prove that Shi arrangements are divisionally free.

\section{Preliminaries}
Let us review several definitions and results used in the rest of 
this article. We use the notation and definitions appeared in \S 1. We use \cite{OT} as a general reference 
in this section.
Define the \textbf{M\"{o}bius function} $\mu:L(\A) \rightarrow \Z$ by 
$
\mu(X)=1$ if $X=V$ and $\mu(X)=-\sum_{X \subsetneq Y \in L(\A)} \mu(Y)$ if 
$X \neq V$. Then the \textbf{characteristic polynomial} $\chi(\A;t)$ and 
the \textbf{Poincar\`{e} polynomial} $\pi(\A;t)$ are defined by 
\begin{eqnarray*}
\chi(\A;t):&=&\sum_{X \in L(\A)} \mu(X) t^{\dim X},\\
\pi(A;t):&=&\sum_{X \in L(\A)} \mu(X) (-t)^{\codim X}.
\end{eqnarray*}
As mentioned in \S1, $\pi(\A;t)$ equals to 
the topological Poincar\`{e} polynomial of $M(\A)$ when 
$\K=\mathbb{C}$. Also, it is easy to see that $\chi(\A)$ is divisible by $(t-1)$ if $\A \neq 
\emptyset$. Write and define 
\begin{eqnarray*}
\chi(\A;t)&=&
\sum_{i=0}^{\ell} b_i(\A)(-1)^i t^{\ell-i},\\
\chi_0(\A;t):&=&\chi(\A;t)/(t-1)=
\sum_{i=0}^{\ell-1} b_i(d\A)(-1)^i t^{\ell-1-i}.
\end{eqnarray*}

\begin{rem}
The reason why we use the terminology $b_i(d\A)$ is as follows.
Fix a hyperplane $H_0 \in \A$. Then we may consider the operation 
named \textbf{deconing} $d\A$ of $\A$ with respect to $H_0$. Namely, 
$d\A$ is a set of affine hyperplanes obtained as intersections of $\alpha_{H_0}=1$ with 
all hyperplanes $L \in \A \setminus \{H_0\}$. Then it is known that 
$
\chi(d\A;t)=\chi_0(\A;t)$. See \cite{OT} for example. This is the reason of the 
notation above.
\end{rem}

Let $\theta_E:=\sum_{i=1}^\ell x_i \partial_{x_i}$ be the Euler derivation. 
Let us recall several results on these polynomials and freeness: 

\begin{theorem}[Deletion-restriction formula, e.g., \cite{OT}, Corollary 2.57]
Let $(\A,\A',\A^H)$ be a triple. Then it holds that 
\begin{eqnarray*}
\chi(\A;t)&=&\chi(\A';t)-\chi(\A^H;t),\\
\chi_0(\A;t)&=&\chi_0(\A';t)-\chi_0(\A^H;t).
\end{eqnarray*}
\label{delres}
\end{theorem}

\begin{theorem}[Addition-deletion theorem, \cite{T1}]
Let $(\A,\A',\A^H)$ be a triple. Then any two of the following three 
imply the third:
\begin{itemize}
\item[(1)]
$\A$ is free with $\exp(\A)=(d_1,\ldots,d_{\ell-1},d_\ell)$.
\item[(2)]
$\A'$ is free with $\exp(\A')=(d_1,\ldots,d_{\ell-1},d_\ell-1)$.
\item[(3)]
$\A^H$ is free with $\exp(\A^H)=(d_1,\ldots,d_{\ell-1})$.
\end{itemize}
\label{adddel}
\end{theorem}

\begin{theorem}[Restriction Theorem, \cite{T1}]
In the notation of Theorem \ref{adddel}, assume that both $\A$ and $\A'$ are 
free. Then all the statements (1), (2), and (3) in Theorem \ref{adddel} hold. Moreover,  
we may choose a basis $\theta_1,\ldots,\theta_\ell$ for $D(\A)$ such that 
$\alpha_H$ divides $\theta_\ell$ and $\pi(\theta_1),\ldots,\pi(\theta_{\ell-1})$ 
form a basis for $D(\A^H)$, where $\pi$ is a residue map $D(\A) \rightarrow 
D(\A^
H)$.
\label{rest}
\end{theorem}

\begin{theorem}[Factorization theorem, \cite{T2}, Main Theorem]
If $\A$ is free with $\exp(\A)=(d_1,\ldots,d_\ell)$, then 
$\chi(\A;t)=\prod_{i=1}^\ell (t-d_i)$.
\label{fact}
\end{theorem}

When $\ell=3$, the addition-deletion theorem \ref{adddel} has the 
following simple form.

\begin{theorem}[\cite{A}, Theorem 1.1 (3)]
Let $\A$ be a central arrangement in $V=\K^3$. 
Then two of the following three imply all the statements (1), (2) and (3) in Theorem \ref{adddel}:
\begin{itemize}
\item[(a)]
$\chi(\A;t)=(t-1)(t-d_1)(t-d_2)$.
\item[(b)]
$\chi(\A';t)=(t-1)(t-d_1)(t-(d_2-1))$.
\item[(c)]
$|\A^H|=d_1+1$.
\end{itemize}
\label{line}
\end{theorem}

Among the equivalent conditions in Theorem \ref{line}, we can find several
divisions of characteristic polynomials, and they imply the freeness of each member in the triple 
$(\A,\A',\A^H)$. From this point of view, Theorem \ref{newadd} is a generalization of 
Theorem \ref{line}. Also, if we regard these freeness as a local freeness in 
codimension three in $\K^3$, then Theorem \ref{codim3} is also a generalization of Theorem \ref{line}. 
%
%
The following is also in \cite{A} with a different formulation.

\begin{theorem}[\cite{A}, Theorem 1.1]
Assume that $\ell=3$ and let us consider the division 
$$
\chi_0(\A;t)=(t-(|\A^H|-1))(t-(|\A|-|\A^H|))+a.
$$
Then $a=\chi_0(\A;|\A^H|-1) \ge 0$. Moreover, 
$a=0$ implies that $\A$ is free.
\label{div}
\end{theorem}

\begin{rem}
Theorem \ref{div} says that the integer $a$ can be regarded as the  
remainder of the polynomial division of $\chi_0(\A;t)$ by $t-(|\A^H|-1)=
\chi_0(\A^H;t)$. 
%
The non-negativity of 
$a$ in the above was 
generalized in \cite{A}, Theorem 7.1 for an arbitrary $\ell \ge 3$ 
just as an inequality. 
However, from the viewpoint of polynomial divisions, 
we may understand the non-negativity as that of  the leading term of 
the remainder of a characteristic polynomial division naturally. 
%
%
Hence the non-negativity $r_0 \ge 0$ in Theorem \ref{codim3} 
can be regarded as a generalization of Theorem \ref{div}. 
\label{r0}
\end{rem}

Now let us explain the freeness criterion by using multiarrangements. A \textbf{multiarrangement} 
$(\A,m)$ 
is a pair consisting of an $\ell$-arrangement $\A$ and a function $m: \A 
\rightarrow \Z_{>0}$. Let $|m|:=
\sum_{
H \in \A} m(H)$. For $L \in \A$, let 
$\delta_L:\A \rightarrow \{0,1\}$ be a \textbf{characteristic multiplicity} of $L$ defined by 
$\delta_L(H)=1$ only when $H=L$ (hence $\delta_L(H)=0$ for all $\A \ni H \neq L)$.
For $(\A,m)$, we can 
define the \textbf{logarithmic derivation module} $D(\A,m)$ by 
$$
D(\A,m):=\{\theta \in \Der S \mid \theta(\alpha_H) \in S \alpha_H^{m(H)}\ 
(\forall H \in \A)\}.
$$
Also, its 
\textbf{freeness} and \textbf{exponents} can be defined in the same manner as for 
$\A$ and $D(\A)$. 
When $(\A,m)$ is free with $\exp(\A,m)=(d_1,\ldots,d_\ell)$, it holds that 
$\sum_{i=1}^\ell d_i=|m|$. For two multiplicities 
$m_i:\A \rightarrow \Z_{>0}\ (i=1,2)$, $m_1 \le m_2$ denotes 
$m_1(H) \le m_2(H)$ for all $H \in \A$. 
Then the following is
the most basic method to determine the freeness. 

\begin{theorem}[Saito's criterion, \cite{Sa}]
Let $\theta_1,\ldots,\theta_\ell \in D(\A,m)$. Then $(\A,m)$ is free if and only if 
$\theta_1,\ldots,\theta_\ell$ are $S$-independent and
$\sum_{i=1}^\ell \deg \theta_i=|m|$.
\label{Saito}
\end{theorem}

Not only the freeness, but also a \textbf{characteristic polynomial} 
$\chi(\A,m;t)$ for a multiarrangement $(\A,m)$ 
can be defined too. Let us write $$
\chi(\A,m;t)=\sum_{i=0}^\ell b_i(\A,m)(-1)^i t^{\ell-i}.
$$
Also, there is the factorization theorem for free multiarrangements as for the 
arrangement cases, i.e., if $(\A,m)$ is free with $\exp(\A,m)=(d_1,\ldots,d_\ell)$, then 
$$
\chi(\A,m;t)=\prod_{i=1}^\ell(t-d_i).
$$
For details, see Definition 2.1 and Theorem 4.1 in \cite{ATW}.
Contrary to $\chi(\A;t)$, to compute $\chi(\A,m;t)$ when $m$ is not identically $1$ 
is very difficult. In fact, it is not combinatorial (see \cite{ATW}). However, 
it is easy to check that $b_1(\A,m)=|m|$. Moreover, $b_2(\A,m)$ can be 
computed by using the following local-global formula.

\begin{theorem}[Local-global formula, \cite{ATW}, Theorem 3.3]
Let $(\A,m)$ be a multiarrangement. For $X \in L_2(\A)$, denote 
$\exp(\A_X,m|_{\A_X})=(d_1^X,d_2^X,0,\ldots,0)$. Then 
$$
b_2(\A,m)=\sum_{X \in L_2(\A)} d_1^Xd_2^X=\sum_{X \in L_2(\A)} b_2(\A_X,m|_{\A_X}).
$$
\label{LG}
\end{theorem}

From a central arrangement $\A$ and a hyperplane $H \in \A$, we may define an induced multiarrangement  $(\A^H,m^H)$ (called the \textbf{Ziegler restriction} of $\A$ onto $H$) 
by 
$$
m^H(K):=|\{L \in \A' \mid L \cap H=K\}|.
$$
The multiplicity $m^H$ is called the \textbf{Ziegler multiplicity}. 
By definition, $|m^H|=|\A'|=|\A|-1$. 
Then Ziegler showed the following fundamental fact.

\begin{theorem}[\cite{Z}]
Assume that $\A$ is free with $\exp(\A)=(1,d_2,\ldots,d_\ell)$. Then 
the Ziegler restriction $(\A^H,m^H)$ onto an arbitrary $H \in \A$ is free with 
$\exp(\A^H,m^H)=(d_2,\ldots,d_\ell)$. 
\label{Ziegler}
\end{theorem}

We say that $\A$ is \textbf{locally free in codimension $i$ along $H \in \A$} if 
$\A_X$ is free for all $X \in L_i(\A)$ with $X \subset H$ 
(equivalently, for all $X \in L_{i-1}(\A^H)$). By using this notion,
as a converse of Theorem \ref{Ziegler}, the following freeness criterion is known.

\begin{theorem}[\cite{AY}, Theorem 5.1]
Let $\A \neq \emptyset$ and $H \in \A$. Then $b_2(d\A) \ge b_2(\A^H,m^H)$, and 
the equality holds if and only if $\A$ is locally free in codimension three along $H$. 
Moreover, 
if the Ziegler restriction $(\A^H,m^H)$ 
is free, then the following conditions are equivalent:
\begin{itemize}
\item[(1)]
$\A$ is free.
\item[(2)]
$b_2(d\A)=b_2(\A^H,m^H)$.
\item[(3)]
$\A_X$ is free for any $X \in L_2(\A^H)$, i.e., 
$\A$ is locally free in codimension three along $H$.
\end{itemize}
\label{criterion}
\end{theorem}

Also, let us recall the addition-deletion theorem for multiarrangements $(\A,m)$ from 
\cite{ATW2}. 
For $H \in \A$, we may define a restriction $(\A^H,m^*)$ 
of the multiarrangement $(\A,m)$ onto $H$, which is called the \textbf{Euler restriction} 
of $(\A,m)$ 
onto $H$. Also, the new multiplicity $m^*$ is called the \textbf{Euler multiplicity}. Let us  
explain the definition of $m^*$ since it will be used for the proofs of main results. 

For $X \in L_2(\A)$ with $X \subset H$, let $m_X$ denote the restriction of $m$ onto $\A_X$. Then 
the multiarrangement $(\A_X,m_X)$ is a direct product of a $2$-multiarrangement 
and the $(\ell-2)$-empty arrangement. It is well-known that every $2$-arrangement is free. Hence 
we may define $\exp(\A_X,m_X)=(d_1^X,d_2^X,0,\ldots,0)$. Therefore, there 
is a homogeneous free basis $\theta_1,\theta_2,\varphi_1,\ldots,\varphi_{\ell-2}$ for $D(\A_X,m_X)$ 
such that $\deg \theta_i=d_i\ (i=1,2)$ and $\deg \varphi_j=0\ (j=1,\ldots,\ell-2)$. By the 
same reason, there are free basis $\theta_1',\theta_2',\varphi_1',\ldots,\varphi_{\ell-2}'$ 
for $D(\A_X,m_X-\delta_H)$ such that
$\deg \theta_i'=d_i'\ (i=1,2)$ and $\deg \varphi_j'=0\ (j=1,\ldots,\ell-2)$, here
$\exp(\A_X,m_X-\delta_H)=(d_1',d_2',0,\ldots,0)$. It was shown in \cite{ATW2} that we may assume 
$d_1=d_1'+1$ and $d_2=d_2'$. Moreover, we may choose bases in such a way that 
$$
\theta_1=\alpha_H \theta_1',\ \theta_2=\theta_2',\ \varphi_j=\varphi_j'\  (j=1,\ldots,\ell-2).
$$
Recall that $X \in L_2(\A)$ with $X \subset H \iff 
X \in \A^H$. Now define 
$m^*(X):=\deg \theta_2$. 
For more details of 
these definitions, see \cite{ATW2}. By using these definitions, 
we have the 
\textbf{addition-deletion theorem for multiarrangements} as follows.

\begin{theorem}[\cite{ATW2}, Theorem 0.8]
Let $(\A,m)$ be a multiarrangement, $H \in \A$ and 
$(\A^H,m^*)$ the Euler restriction of $(\A,m)$ onto $H$. Then any two of the following three 
imply the third:
\begin{itemize}
\item[(1)]
$(\A,m)$ is free with $\exp(\A,m)=(d_1,\ldots,d_{\ell-1},d_\ell)$.
\item[(2)]
$(\A,m-\delta_H)$ is free with $\exp(\A,m-\delta_H)=(d_1,\ldots,d_{\ell-1},d_\ell-1)$.
\item[(3)]
$(\A^H,m^*)$ is free with $\exp(\A^H,m^*)=(d_1,\ldots,d_{\ell-1})$.
\end{itemize}
Moreover, if both $(\A,m)$ and $(\A,m-\delta_H)$ are free, then all the 
three hold true.
\label{multiadd}
\end{theorem}

\section{Proofs of results}
Let us start the proof of main results in \S 1. Let us add one notation.
For an $\ell$-multiarrangement $(\A,m)$ and $X \in L_i(\A)$, assume that 
$(\A_X,m_X)$ is free with $\exp(\A_X,m_X)=(d_1,\ldots,d_i,0,\ldots,0)$. Then the set 
$\exp^*(\A_X,m_X):=\{d_1,\ldots,d_i\}$. 

First,  we introduce a lemma which is used without referring in the 
rest of this article.

\begin{lemma}
A $2$-multiarrangement $(\A,m)$ is free. 
Let $\exp(\A,m)=(d_1,d_2)$. Then 
for $H \in \A$, it holds that 
$\exp(\A,m-\delta_H)=(d_1-1,d_2)$ or $(d_1,d_2-1)$.
\label{basic}
\end{lemma}

\noindent
\textbf{Proof}. For example, see Lemma 2.7 in \cite{A}.\owari
\medskip

Second we need 
the following lemmas. 

\begin{lemma}
Let $\A$ be a $2$-arrangement. Then 
$\exp(\A,m)=(|m|-|\A|+1,|\A|-1)$ when 
$|m| \le 2|\A|-1$. If $|m|  \ge 2|\A|$, then $d_i \ge |\A|-1\ (i=1,2)$ 
for $\exp(\A,m)=(d_1,d_2) $ with $d_1 \le d_2$.
\label{2free}
\end{lemma}

\noindent
\textbf{Proof}. 
This result is well-known. Here we give a 
short proof. 
For a complete proof, see 
Lemma 2.10 in \cite{A} for example.
%

First consider the case $|m| \le 2|\A|-1$. 
Then clearly $$
\theta_1:=(\prod_{H \in \A} \alpha_H^{m(H)-1}) \theta_E \in D(\A,m).
$$
Note that $\deg \theta_1=|m|-|\A|+1$. 
Also the definition of $D(\A,m)$ shows that there are no $\theta \in D(\A,m)$ such that 
$\theta \mid \theta_1$. Since $d_1+d_2=|m|$ by Theorem \ref{Saito}, 
it has to hold that $|m|-|\A|+1=d_1$, or $|m|-|\A|+1=d_2$, or  
$d_2 <|m|-|\A|+1$. The third case  contradicts Saito's criterion \ref{Saito}.

Second assume that $|m| \ge 2|\A|$ and take multiplicities $m' \le m$ such that 
$|m'| =2|\A|-1$. Then the first assertion shows that 
$\exp(\A,m')=(|\A|-1,|\A|)$. Since 
$D(\A,m') \supset D(\A,m)$, it holds that 
$|\A|-1 \le d_1 \le d_2$.\owari
\medskip

\begin{lemma}
Let $(\A,m)$ be a multiarrangement and fix $H \in \A$ with $m(H) \ge 2$. 
Then 

(1)\,\,
if $\A$ is a $2$-arrangement, then $b_2(\A,m)-b_2(\A,m-\delta_H) =
\exp^*(\A,m) \cap \exp^*(\A,m-\delta_H)$, and 

(2)\,\,
$b_2(\A,m)-b_2(\A,m-\delta_H) \ge |\A|-1$.
\label{bound}
\end{lemma}

\noindent
\textbf{Proof}.  

(1)\,\,
Immediate by the definition of $b_2(\A,m)$ and 
Lemmas \ref{basic} and \ref{2free}.

(2)\,\,
Let us compute $b_2(\A,m)-b_2(\A,m-\delta_H)$ by using the local-global 
formula \ref{LG}. Since $(\A_X,m_X)=(\A_X,(m-\delta_H)_X)$ for $L_2(\A) \ni X \not \subset H$,
it holds that
\begin{equation}
b_2(\A,m)-b_2(\A,m-\delta_H)=\sum_{X \in \A^H} (b_2(\A_X,m_X)-b_2(\A_X,(m-\delta_H)_X)).
\label{eq101}
\end{equation}
Now apply (1) and Lemma \ref{2free} to show that 
$$
(\ref{eq101})
\ge \sum_{X \in \A^H} (|\A_X|-1)=
|\A|-1,
$$
which completes the proof.\owari
\medskip

\medskip

\begin{lemma}
Let $(\A,m)$ be a multiarrangement and fix $H \in \A$ with $m(H) \ge 2$. Let 
$m^*$ be the Euler multiplicity of $(\A,m)$ onto $H$. Then
 
(1)\,\,
$m^*(X) \ge |\A_X|-1$ for all $X \in \A^H$, 
$m^*(X)=\exp^*(\A_X,m_X)\cap \exp^*(\A_X,m_X-\delta_H)$ and 

(2)\,\,
if 
$m^*(X) \neq |\A_X|-1$ for some $X \in \A^H$, then 
$b_2(\A,m)-b_2(\A,m-\delta_H)>|\A|-1$.
\label{strict}
\end{lemma}

\noindent
\textbf{Proof}. 
(1)\,\,
When $|m_X| \le 2|\A_X|-1$, the explicit construction of $\theta_1$ in 
the proof of Lemma \ref{2free} shows that $m^*(X)=|\A_X|-1$. If 
$|m_X| \ge 2|\A_X|$, then 
Lemma \ref{2free} shows that 
$d \ge |\A_X|-1$ for all $d \in \exp^*(\A_X,m_X)$. Hence the definition of 
$m^*$ shows that $m^*(X) \ge |\A_X|-1$. The second statement follows immediately by 
the same proof as Lemma \ref{bound} (1).


(2)\,\,
By the assumption, 
(1) and Lemma \ref{bound} (1), 
$m^*(Y)=b_2(\A_Y,m_Y)-b_2(\A_Y,(m-\delta_H)_Y)>|\A_Y|-1$ for some $Y \in \A^H$. 
Then apply this and Lemma \ref{bound} (2) to the local-global formula \ref{LG} to obtain 
\begin{eqnarray*}
b_2(\A,m)-b_2(\A,m-\delta_H)&=&\sum_{X \in \A^H}( b_2(\A_X,m_X)
-b_2(\A_X,(m-\delta_H)_X))\\
&>& \sum_{X \in \A^H} (|\A_X|-1)=|\A|-1.
\end{eqnarray*}
\owari\medskip
%

\begin{lemma}
Let $(\A,m)$ be an $\ell$-multiarrangement, $H \in \A$ with $m(H)>1$, 
$X \in L_2(\A^H)$ and $(\A^H,m^*)$ the Euler restriction of $\A$ onto $H$. Then 
\begin{equation}
b_2(\A_X^H,m_X)-b_2(\A^H_X,(m-\delta_H)_X)=m^*(X).
\end{equation}
\label{eq}
\end{lemma}

\noindent
\textbf{Proof}. Immediate from Lemmas \ref{bound} (1) and 
\ref{strict} (1).
\owari
\medskip

\begin{prop}
In the setup of Theorem \ref{codim3}, it holds that 
\begin{eqnarray}
&\ &b_2(\A^H)+(|\A^H|-1)(|\A|-|\A^H|-1)+r_0 \\
&=&b_2(d\A) \ge b_2(\A^H,m^H) \nonumber \\
&\ge &b_2(\A^H)+(|\A^H|-1)(|\A|-|\A^H|-1). \nonumber
\label{ineq}
\end{eqnarray}
In particular, when $r_0=0$, for a multiplicity $m$ on $\A^H$ such that 
$1 \le m \le m^H$,
it holds that 
\begin{eqnarray}
\ \ \ \ b_2(\A^H,m)-b_2(\A^H,m-\delta_L)&=&|\A^H|-1,\label{1}\\
\ \ \ \ b_2(\A_X^H,m_X)-b_2(\A_X^H,(m-\delta_L)_X
)&=&\label{2}|\A^H_X|-1
\end{eqnarray}
for $L \in \A^H$, 
$X \in L_2(\A^H)$ with $X \subset L$ and 
$m(L)>1$.
\label{eqeq}
\end{prop}

\noindent
\textbf{Proof}.
Recall the division of characteristic polynomials:
$$
\chi_0(\A;t)=\chi_0(\A^H;t)(t-(|\A|-|\A^H|))+\sum_{i=0}^{\ell-3} r_i t^{\ell-3-i}.
$$
Now let us 
introduce three (in)equalities:

\textbf{First inequality}. $b_2(d\A)\ge b_2(\A^H,m^H)$. This is a part of Theorem \ref{criterion}. 

\textbf{Second inequality}.  $b_2(\A^H,m^H) \ge b_2(\A^H)+(|\A^H|-1)(|\A|-|\A^H|-1)$. 
First note that $|m_H|-|\A^H|=|\A|-|\A^H|-1$. By definition, there are 
a sequence of hyperplanes $L_1,\ldots,L_{|\A|-|\A^H|-1} \in \A^H$ such that 
$m_H=1+\delta_{L_1}+\cdots+\delta_{L_{|\A|-|\A^H|-1}}$. 
Hence Lemma \ref{bound} (2) implies that
\begin{eqnarray}
\label{key}b_2(\A^H,m^H)&=&
\sum_{i=1}^{|\A|-|\A^H|-1}(b_2(\A^H,
1+\delta_{L_1}+\cdots+\delta_{L_i})\\
&\ &-b_2(\A^H,1+\delta_{L_1}+
\cdots+\delta_{L_{i-1}})) 
+b_2(\A^H) \nonumber \\
&\ge&\sum_{i=1}^{|\A|-|\A^H|-1}(|\A^H|-1)+b_2(\A^H)\nonumber\\
&=& b_2(\A^H)+(|\A|-|\A^H|-1)(|\A^H|-1). \nonumber 
\end{eqnarray}





\textbf{Third equality}. 
$b_2(d\A)=b_2(\A^H)+(|\A^H|-1)(|\A|-|\A^H|-1)+r_0$. First, note that 
$\chi_0(\A;t)=(t-(|\A|-|\A^H|))\chi_0(\A^H;t)+\sum_{i=0}^{\ell-3}r_it^{\ell-3}$. 
Hence comparing the coefficients of $t^{\ell-3}$ implies 
$$
b_2(d\A)=(|\A|-|\A^H|)(|\A^H|-1)+b_2(d\A^H)+r_0.
$$
Since $\chi(\A^H;t)=(t-1)\chi_0(\A^H;t)$, it holds that 
$b_2(\A^H)=b_2(d\A^H)+|\A^H|-1$. 
Hence 
$$
b_2(d\A)=(|\A|-|\A^H|-1)(|\A^H|-1)+b_2(\A^H)+r_0.
$$
Now we have three (in)equalities. Combine these three to obtain 
\begin{eqnarray*}
&\ &b_2(\A^H)+(|\A^H|-1)(|\A|-|\A^H|-1)+r_0 \\
&=&b_2(d\A) \ge b_2(\A^H,m^H) \\
&\ge &b_2(\A^H)+(|\A^H|-1)(|\A|-|\A^H|-1),
\end{eqnarray*}
which is the first statement in this proposition.

Now assume that $r_0=0$. Then the inequalities above are all equalities. Hence 
\begin{eqnarray}
b_2(d\A)&=&b_2(\A^H,m^H) \label{eq7}\\
&=&b_2(\A^H)+(|\A|-|\A^H|-1)(|\A^H|-1).
\label{eq6}
\end{eqnarray}

Let us show
$b_2(\A^H,m)-
b_2(\A^H,m-\delta_L)=|\A^H|-1$ for a multiplicity $1 \le m \le m^H$ and 
$L \in \A^H$ with $m(L)>1$.
We know that 
$b_2(\A^H,m)-
b_2(\A^H,m-\delta_L)\ge |\A^H|-1$ by Lemma \ref{bound} (2). If this 
inequality is strict, then 
again Lemma \ref{bound} (2) and (\ref{key}) show that the equality (\ref{eq6}) cannot 
hold, which is a contradiction.

Next let us show that 
$b_2(\A^H_X,m_X)-
b_2(\A^H_X,(m-\delta_L)_X)=|\A^H_X|-1$ for $X \in L_2(\A^H)$ with 
$X \subset L$. By the above, we know that
$b_2(\A^H,m)-
b_2(\A^H,m-\delta_L)=|\A^H|-1$. So the local-global formula \ref{LG} shows that 
\begin{eqnarray*}
b_2(\A^H,m)-
b_2(\A^H,m-\delta_L)&=&
\sum_{X \in (\A^H)^L} (b_2(\A_X^H,m_X)-
b_2(\A_X^H,(m-\delta_L)_X))\\
&=&|\A^H|-1=
\sum_{X \in (\A^H)^L}  (|\A_X^H|-1).
\end{eqnarray*}
Since $b_2(\A_X^H,m_X)-
b_2(\A_X^H,(m_X-\delta_L)_X)\ge |\A_X^H|-1$ by Lemma \ref{bound} (2), 
this inequality has to be an equality, which completes the proof. \owari
\medskip

Now let us prove Theorem \ref{codim3}. 
\medskip

\noindent
\textbf{Proof of Theorem \ref{codim3}}. 
The nonnegativity $r_0 \ge 0$ holds immediately from 
the inequality in Proposition \ref{eqeq}.
Assume that $r_0=0$. Then we have the 
equation (\ref{eq7}). 
Thus Theorem \ref{criterion} shows that $\A$ is locally free in 
codimension three along $H$. Hence to complete the proof by applying Theorem \ref{criterion}, 
it suffices to show that 
$\chi(\A^H_X;t)$ divides $\chi(\A_X;t)$ for all $X \in L_3(\A)$ with $X \subset H$. 

We show that $\exp(\A_X^H,(m^H)_X)=
(|\A_X^H|-1,|(m^H)_X|-|\A^H_X|+1)$ by applying Theorem \ref{multiadd}. 
First, note that 
$((\A^H)_X,(m^H)_X)$ coincides with 
the Ziegler restriction of $\A_X$ onto $H \supset X$ by definition. Hence let us express this 
multiarrangement by $(\A_X^H,m_X^H)$. 
Now we need to compute 
the Euler multiplicity $m^*$, which can be obtained as the Euler 
restriction of $\A_X$ onto $H$. 
Recall that we have the equation (\ref{2}). Hence Lemma \ref{eq} shows that
\begin{equation}
m^*(X)=b_2(\A_X^H,m_X)-b_2(\A_X^H,(m-\delta_L)_X)=|\A_X^H|-1
\label{eq10}
\end{equation} 
for $L \in \A^H$ and $X \in (\A^H)^L$ as in the notation of Proposition \ref{eqeq}.

Since $\exp^*(\A^H_X)=(1,|\A_X^H|-1)$ by 
Lemma \ref{2free}, Theorem \ref{multiadd} 
shows that 
$\exp(\A_X^H,(m^H)_X)=(|(m^H)_X|-(|\A_X^H|-1), |\A_X^H|-1)
=(|\A_X|-|\A_X^H|, |\A_X^H|-1)$. 
Because $\A$ is locally free in codimension three along $H$ by (\ref{eq7}),  
Theorems \ref{fact} and \ref{Ziegler} show that 
\begin{eqnarray*}
\chi(\A_X;t)&=&(t-1)(t-(|\A_X^H|-1))(t-(|\A_X|-|
\A_X^H|)),\\ 
\chi(\A_X^H;t)&=&(t-1)(t-(|\A_X^H|-1)),
\end{eqnarray*}
which completes the proof. \owari
\medskip

\begin{example}
Let $\A$ be a plane arrangement consisting of the cone of 
all edges and diagonal lines of a regular pentagon. 
Hence $|\A|=11$. It is well-known that (for example, see \cite{OT}, 
Example 4.54)  
$\A$ is free with $\exp(\A)=(1,5,5)$ and $|\A^H|=5$ for any $H \in \A$. 
Let $\B$ be the coning of $\A$, hence free with $\exp(\A)=(1,1,5,5)$. 
Let $H_0 \in \B$ be the infinite line in $\A$ and let $X \in L_2(\A^{H_0})$ 
be a flat which is contained in all the cone of planes belonging to $\A$. Then it 
is easy to see that
$\B_X=\A \times \emptyset_1$ and 
$\B_X^{H_0}$ consists of five linear lines in a plane. Hence 
$$
\chi(\B_X;t)=(t-1)(t-5)^2,\ 
\chi(\B_X^{H_0};t)=(t-1)(t-4).
$$
So clearly $\chi(\B_X^{H_0};t) \nmid \chi(\B_X;t)$.
%
\label{nondiv}
\end{example}



\noindent
\textbf{Proof of Theorem \ref{newadd}}. 
By Theorems \ref{codim3} and \ref{criterion}, it suffices to show that the 
Ziegler restriction $(\A^H,m^H)$ of $\A$ onto $H$ is free. By the assumption and 
Terao's factorization theorem \ref{fact}, it holds that 
$\chi(\A^H;t)=\prod_{i=1}^{\ell-1}(t-d_i)$, where 
$d_1=1$ and $\exp(\A^H)=(d_1,\ldots,d_{\ell-1})$. Since $\chi(\A^H;t) \mid
\chi(\A;t)$, we may write 
$\chi(\A;t)=(t-d_\ell)\chi(\A^H;t)=\prod_{i=1}^\ell (t-d_i)$ for $d_\ell:=|\A|-|\A^H| \in \Z_{\ge 0}$. 
In fact, we show that, 
for every multiplicity $m$ with $1 \le m \le m^H$, the multiarrangement $(\A^H,m)$ is free with 
$\exp(\A^H,m)=(|m|-|\A^H|+1,d_2,\ldots,d_{\ell-1})$.

We show by induction on $|m|$. If $m$ is a constant multiplicity $1$, then 
the statement holds true by the assumption. Assume that the statement holds true 
for all $m$ with $|m|-|\A^H|<n\ ( n \in \Z_{\ge 1})$. Take an arbitrary $m$ 
satisfying the condition above and $|m|-|\A^H|=n$.  
We may choose a hyperplane $L \in \A^H$ such that $m(L) \ge 2$. 
By the induction hypothesis, the multiarrangement $(\A^H,m-\delta_L)$ is free with 
$\exp(\A^H,m-\delta_L)=(n,d_2,\ldots,d_{\ell-1})$. Let us apply 
the addition theorem for multiarrangements (Theorem \ref{multiadd}) 
to confirm that $(\A^H,m)$ is free with $\exp(\A^H,m)=(n+1,d_2,\ldots,d_{\ell-1})$.

For that purpose, we need to compute the Euler restriction $((\A^H)^L,m^*)$
of $(\A^H,m)$ onto $L$. 
By the assumption, Proposition \ref{eqeq}, the equation (\ref{2}) and 
Lemma \ref{eq}, 
the equation (\ref{eq10}) also follows here, i.e., 
$$
b_2(\A^H_X,m_X)-b_2(\A^H_X,(m-\delta_L)_X)=m^*(X)=|\A_X^H|-1
$$
for all $X \in (\A^H)^L$.
%
On the other hand, by definition, $m^L(X)=|\A^H_X|-1$ for 
$X \in (\A^H)^L$, where $((\A^H)^L,m^L)$ is 
the Ziegler restriction of $\A^H$ onto $L$. Hence in this case, 
these two restrictions coincide; $((\A^H)^L,m^L)=((\A^H)^L,m^*)$. 
Since 
$\A^H$ is free with $\exp(\A^H)=(1,d_2,\ldots,d_{\ell-1})$, Theorem \ref{Ziegler} 
shows that $((\A^H)^L,m^L)=((\A^H)^L,m^*)$ is free with exponents 
$(d_2,\ldots,d_{\ell-1})$. Hence the assumption on the freeness and 
exponents on $(\A^H,m-\delta_L)$, and Theorem \ref{multiadd} show that 
$(\A^H,m)$ is also free with $\exp(\A^H,m)=(n+1,d_2,\ldots,d_{\ell-1})$. 

Now apply this argument to $(\A^H,m^H)$. Then we obtain that $(\A^H,m^H)$ is 
free with 
\begin{eqnarray*}
\exp(\A^H,m^H)&=&(|m^H|-|\A^H|+1,d_2,\ldots,d_{\ell-1})\\
&=&(|\A|-|\A^H|,d_2,\ldots,d_{\ell-1})\\
&=&(d_2,\ldots,d_{\ell-1},d_\ell), 
\end{eqnarray*}
which completes the proof. When 
the condition $\chi(\A^H,t) \mid \chi(\A;t)$ is replaced by 
$r_0=0$ as mentioned in Remark \ref{r0division}, the proof is totally the same, so we left it to the reader. 
\owari
\medskip

\noindent
\textbf{Proof of Theorem \ref{tower}}. Apply
Theorem \ref{newadd} repeatedly with Lemma \ref{basic}. \owari
\medskip


\begin{rem}
Theorem \ref{newadd} can be regarded as a modified converse of the famous 
Orlik's conjecture, which asserted that $\A^H$ is free if $\A$ is free. 
To this conjecture, a counter example was found by Edelman and Reiner in \cite{ER}. 
This conjecture was asserting that global freeness implies restricted freeness. Though this is 
not true, Theorem \ref{newadd} asserts that the modified 
converse is true, i.e., restricted freeness with 
a combinatorial condition implies global freeness. 
\label{Orlik}
\end{rem}

As a corollary of the proof above, the following is immediate.

\begin{cor}
(1)\,\,
Assume that $\chi(\A^H;t) \mid \chi(\A;t)$. Then $(\A^H,m^H)$ is free if 
$\A^H$ is free.
\newline
(2)\,\,
Assume that there is $H \in \A$ such that 
the Ziegler restriction $(\A^H,m^H)$ is free and 
$\chi(\A^H;t) \mid \chi(\A;t)$. Then $\A$ is free.
\newline
(3)\,\,
Assume that $\A^H$ is free. If $r_0=0$ in the division 
$$
\chi_0(\A;t)=(t-(|\A|-|\A^H|))\chi_0(\A^H;t)+\sum_{i=0}^{\ell-3} r_i t^{\ell-3-i},
$$
then $r_1=\cdots=r_{\ell-3}=0$. Equivalently, $\chi(\A^H;t)$ divides 
$\chi(\A;t)$.
\label{simplemulti}
\end{cor}


If we emphasize the aspect as the addition-deletion theorems, we may have the following 
formulations too.

\begin{theorem}
Let $\A$ be an $\ell$-central arrangement. Assume that there is a 
hyperplane $L \not \in \A$ such that $\chi((\A \cup \{L\})^L;t)$ divides $\chi(\A;t)$ and that 
$(\A \cup \{L\})^L$ is free. Then 
$\A$ is free.
\label{adddiv}
\end{theorem}

\noindent
\textbf{Proof}. Apply the same proof as that of Theorem \ref{newadd} to $\A \cup \{L\}$ to 
obtain that $\A \cup \{L\}$ is free since 
$\chi((\A \cup \{L\})^L;t) \mid \chi( \A \cup \{L\};t)$ by Theorem \ref{delres}. Now apply 
Theorem \ref{adddel} to complete the proof. \owari
\medskip

\begin{theorem}
Let $(\A,\A',\A^H)$ be the triple with respect to $H \in \A$. Assume that 
$\A^H$ is free. Then the following conditions are equivalent:
\begin{itemize}
\item[(1)]
$\A$ is free and $\exp(\A) \supset \exp(\A^H)$.
\item[(2)]
$\A'$ is free and $\exp(\A') \supset \exp(\A^H)$.
\item[(3)]
All the three of $\A,\A'$ and $\A^H$ are free.
\item[(4)]
$\chi(\A^H;t)$ divides $\chi(\A;t)$.
\item[(5)]
$\chi(\A^H;t)$ divides $\chi(\A';t)$.
\item[(6)]
$\chi(\A;t)$ and $\chi(\A';t)$ have a GCD of degree $\ell-1$.
\item[(7)]
In the division 
$\chi_0(\A;t)=\chi_0(\A^H;t)(t-(|\A|-|\A^H|))+\sum_{i=0}^{\ell-3} r_i t^{\ell-3-i}$, 
it holds that $r_0=0$. 
\item[(8)]
In the division 
$\chi_0(\A';t)=\chi_0(\A^H;t)(t-(|\A'|-|\A^H|))+\sum_{i=0}^{\ell-3} r_i' t^{\ell-3-i}$, 
it holds that $r_0'=0$. 
\end{itemize}
\label{same}
\end{theorem}

\noindent
\textbf{Proof}. The equivalences among (1), (2) and (3) follow immediately from 
the addition-deletion theorem \ref{adddel}. Those among (4), (5) and (6) follow from 
the deletion-restriction theorem \ref{delres}. Since $\chi_0(\A;t)=\chi_0(\A';t)-\chi_0(\A^H;t)$, 
that between (7) and (8) is easy. Both 
$(1) \iff (4)$ and $(1)\iff (7)$ follow by Theorem \ref{newadd} and Remark \ref{r0division}, which 
completes the proof. \owari
\medskip

The conditions (1), (2) and (3) are the addition-deletion theorems that 
contain freeness conditions. However, the others are just combinatorial ones. Also, 
the freeness is assumed for an $(\ell-1)$-arrangement $\A^H$ to check 
the freeness of $\ell$-arrangement $\A$, which enables us an inductive argument.


The following result is a corollary of the proofs of Theorems \ref{codim3} 
and \ref{newadd}.


\begin{cor}
Let $\A$ be a free $\ell$-arrangement with $\exp(\A)=(1,d_2,\ldots,d_\ell)$ 
and let $m:\A \rightarrow 
\Z_{>0}$. Let us consider the division 
$$
\chi(\A,m;t)=(t-(|m|-|\A|+1))\chi_0(\A;t)+\sum_{i=0}^{\ell-2} s_it^{\ell-2-i}.
$$
If $s_0=0$, then $(\A,m)$ is free with $\exp(\A,m)=
(|m|-|\A|+1,d_2,\ldots,d_\ell)$.
\label{multiaddfree}
\end{cor}

The following is a certain converse of Theorem \ref{codim3}.

\begin{cor}
Let $\A$ be an $\ell$-arrangement.
\begin{itemize}
\item[(1)]
Assume that there is $H \in \A$ such that 
$\chi(\A^H_X;t)$ divides $\chi(\A_X;t)$ for all $X \in L_2(\A^H)$, and 
$\A^H$ is free. Then $\A$ is free.
\item[(2)]
If there is a flag
$V=X_0 \supset X_1 \supset \cdots \supset X_{\ell-2}
$
with $X_i \in L_i(\A)\ (i=0,\ldots,\ell-2)$ 
such that, for $\A_{i}:=\A^{X_i}\ (i=0,\ldots,\ell-2)$, 
$\chi((\A_{i+1})_Y;t)$ dividies 
$\chi((\A_{i})_Y;t)$ for all $Y \in L_2(\A_{i+1})$ and for $i=0,\ldots,\ell-3$. 
Then $\A$ is free.
\end{itemize} 
\label{tower2}
\end{cor}

\noindent
\textbf{Proof}.
(1)\,\,
By the same argument as in the proof of Theorem \ref{newadd}, it holds that 
$$
b_2(\A,m)-b_2(\A,m-\delta_L)=|\A^H|-1
$$
and
$$
b_2(\A^H_X,m_X)-b_2(\A^H_X,(m-\delta_L)_X)=m^*(X)=|\A_X^H|-1
$$
for all $m:\A^H \rightarrow \Z$ with $1 \le m \le m^H$, 
all $L \in \A^H$ with $m(L)>1$ and $X \in (\A^H)^L$. Hence 
Theorems \ref{Ziegler} and \ref{multiadd} complete the proof. (2) is an 
immediate consequence of (1). \owari
\medskip

Without the freeness of $\A^H$, Corollary \ref{tower2} (1) does not hold. See 
the following example.

\begin{example}
Let $\A:=\{xyzw(x+y+z+w)=0\}$ and $H:=\{w=0\} \in \A$. Then it is easy to check that 
$\chi(\A_X;t)=(t-1)^3$ and $\chi(\A_X^H;t)=(t-1)^2$ for all 
$X \in L_2(\A^H)$. Also, it is easy to 
compute 
$$
\chi(\A;t)=(t-1)(t^3-4t^2+6t-4),\ \ \ 
\chi(\A^H;t)=(t-1)(t^2-3t+3).
$$
Hence $\A^H$ is not free and $\chi(\A^H;t) \nmid \chi(\A;t)$, and Theorem \ref{tower2} (1) does not hold. 
\end{example}

\section{Divisionally free arrangements}

In this section we give the main application of main results in this article.
Before that, let us recall Terao's conjecture. 

\begin{conj}[Terao]
Let $\A,\B$ be two $\ell$-arrangements such that 
$L(\A) \simeq L(\B)$ as posets and that $\A$ is free. Then $\B$ is also free.
\label{TC}
\end{conj}

Conjecture \ref{TC} is still open even when $\ell=3$. 
To approach Conjecture \ref{TC}, it is one way to construct a class of free 
arrangements such that whether an arrangement belongs to that class or not 
can be determined by combinatorics. If we can find such a class and show that all 
free arrangements are in it, then Conjecture \ref{TC} is settled positively. Also, 
how many free arrangements could belong to such a class is an interesting problem itself. 
For that purpose, 
Terao introduced 
a class of inductively free arrangements by using Theorem \ref{adddel} 
which satisfies the properties written above. 
First recall the definition of the inductively free arrangements.

\begin{define}[\cite{T1}]
The set $\mathcal{IF}_\ell$ of $\ell$-arrangements consists of the following:
\begin{itemize}
\item[(1)] If $\ell=1,2$, then all arrangements are in $\mathcal{IF}_\ell$.
\item[(2)] Assume that $\ell \ge 3$. Then 
$\A \in \mathcal{IF}_\ell$ if $\A=\emptyset_\ell$, or 
there is $H \in \A$ such that 
$\A':=\A \setminus \{H\} \in \mathcal{IF}_\ell$,
$\A^H \in \mathcal{IF}_{\ell-1}$ and $\chi(\A^H;t) \mid \chi(\A';t)$.
\end{itemize}
By Theorem \ref{adddel}, $\A \in \mathcal{IF}_\ell$ is free. 
An arrangement $\A$ is called \textbf{inductively free (IF)} if $\A \in 
\mathcal{IF}:=\cup_{\ell \ge 1} \mathcal{IF}_\ell$.
\label{IF}
\end{define}

Hence to approach Conjecture \ref{TC}, 
to enlarge the class $\mathcal{IF}$ is an interesting strategy. 
For that purpose, we introduce the class of divisionally free arrangements.
Though it 
has been already defined in Definition \ref{DF0}, let us give an 
another equivalent definition of it similar to Definition \ref{IF}. 

\begin{definetheorem}
The set $\mathcal{DF}_\ell$ of $\ell$-arrangements can be also defined as follows:
\begin{itemize}
\item[(1)] When $\ell=1$ and $2$, all arrangements are in $\mathcal{DF}_\ell$.
\item[(2)] Assume that $\ell \ge 3$. Then 
$\A \in \mathcal{DF}_\ell$ if $\A=\emptyset_\ell$, or 
there is $H \in \A$ such that 
$\chi(\A^H;t) \mid \chi(\A;t)$ and 
$\A^H \in \mathcal{DF}_{\ell-1}$.
\end{itemize}
An arrangement $\A$ is called \textbf{divisionally free (DF)} if $\A \in 
\mathcal{DF}:=\cup_{\ell \ge 1} \mathcal{DF}_\ell$.
\label{DF}
\end{definetheorem}

Comparing Definition \ref{IF} with Definition \ref{DF}, it can be seen that 
the divisionally free arrangement is an analogy of the inductively free arrangement 
by replacing the addition-deletion theorem \ref{adddel} by Theorem \ref{newadd}.  


By the Definitions \ref{DF0}, \ref{DF}, Theorems \ref{adddel}, \ref{newadd} and 
\ref{tower}, 
the following is clear.

\begin{theorem}
(1)\,\,
If $\A \in \mathcal{DF}$, then $\A$ is free.\\
(2)\,\,
$\mathcal{IF} \subset \mathcal{DF}$.\\
(3)\,\,
Whether $\A \in \mathcal{DF}$ or not depends only on $L(\A)$. 
\label{DFTC}
\end{theorem}

Hence to check the divisional freeness of an arrangement gives one method to 
check the ``dependency of that freeness on combinatorics''.

Recall that every inductively free arrangement has to be constructed from 
the empty arrangement by the addition theorem. On the other hand, 
a divisionally free arrangement need not. 
In fact, the class $\mathcal{DF}$ is strictly larger than $\mathcal{IF}$. 
Let us show it as Theorem \ref{IFDF}.
\medskip

\noindent
\textbf{Proof of Theorem \ref{IFDF}}.
It suffices to find an arrangement $\A \in \mathcal{DF} \setminus \mathcal{IF}$. 
We can find such an example among 
the reflecting hyperplanes of a unitary reflection group. 

Let 
$G_{31}$ be a finite unitary reflection group acting on $\mathbb{C}^4$, where we use 
the labeling of such groups due to Shephard and Todd in \cite{ST}. Let 
$\A$ be the unitary reflection arrangement in $V=\mathbb{C}^4$ corresponding to 
$G_{31}$. Then it is shown that $\A$ is free with $\exp(\A)=(1,13,17,29)$ (see 
\cite{OT}, Table C. 12 for example), but $\A$ is not inductively free 
(\cite{HR}, Theorem 1.1). 

However, Lemma 3.5 in \cite{HR}
showed that there is $H \in \A$ such that $\A^H$ is free with 
$\exp(\A^H)=(1,13,17)$. Also, Lemma 4.1 in \cite{AHR} shows that 
there is $L \in \A^H$ such that $|(\A^H)^L|=14$. 
Hence Theorem \ref{fact} shows that $\chi(\A^H;t) 
\mid \chi(\A;t)$ and 
$\chi((\A^H)^L;t) 
\mid \chi(\A^H;t)$, which implies that $\A$ is divisionally free by 
Theorem \ref{tower}. Therefore, 
$\mathcal{IF} \subsetneq \mathcal{DF}$. 

Also, the intermediate arrangements $\A_\ell^k(r) \in \mathcal{DF} \setminus \mathcal{IF}$ if 
$k \neq 0,\ \ell \ge 3$. See Theorem \ref{akr} for details.
\owari
\medskip


One example to which the divisional freeness is applied is as follows:

\begin{example}
It is known that all the freeness of Weyl arrangements depends only on its combinatorics.
Here we show the same statement, for certain cases, by using 
divisionally free arrangements.

Let $\B_\ell$ be the Weyl arrangement of the type $B_\ell$ defined by 
$$
\prod_{i=1}^\ell x_i \prod_{1 \le i < j\le \ell} (x_i^2-x_j^2)=0.
$$
Then immediately $\B_\ell^{x_\ell=0} =\B_{\ell-1}$. Also, 
we may compute 
$$
\chi(\B_\ell;t)=(t-1)(t-3)\cdots (t-(2\ell-1)).
$$
Hence 
it follows that 
$$
\chi(\B_{i-1};t) \mid \chi(\B_i;t)\ (i=3,\ldots,\ell).
$$
Therefore, $\B_\ell \in \mathcal{DF}$.
The same proof works for the root systems of the types $A$ and $C$.
\end{example}

\begin{example}
By the proof of Theorem \ref{IFDF}, $\A(G_{31}) \in \mathcal{DF}$. 
\label{31}
\end{example}

For the finite unitary reflection arrangements, 
inductive freeness and divisional freeness almost coincide:

\begin{cor}
Let $W$ be a finite irreducible complex reflection group and 
$\A=\A(W)$ its corresponding reflection arrangement. Then 
$\A \in \mathcal{IF}$ or $W=G_{31}$ if and only if 
$\A \in \mathcal{DF}$.
\label{312}
\end{cor}

\noindent
\textbf{Proof}. 
Let $W \neq G_{31}$. Then Corollary 1.3 in \cite{HR} shows that $\A=\A(W) \in \mathcal{IF}$ 
if and only if $\exp(\A^H) \subset \exp(\A)$ for any $H \in \A$. Since $\A$ is free, 
this is 
also equivalent to say that $\A \in \mathcal{DF}$. Hence the proof of 
Theorem \ref{IFDF} completes the proof. \owari
\medskip

It seems interesting to study the (non-)divisional freeness of 
several non-inductively free arrangements around complex reflection arrangements which 
appeared in the recent development. See \cite{AHR}, \cite{HR} for example. 

Recalling the fact that there have been several developments on 
how to check the inductively freeness of several arrangements due to 
\cite{AHR}, \cite{HR} and so on, practically, 
the following formulation of Theorem \ref{tower} is also 
useful.

\begin{theorem}
$\A \in \mathcal{DF}$ if and only if 
$\A$ has a flag
$
V=X_0 \supset X_1 \supset \cdots \supset X_{\ell-2}$ with $X_i \in 
L_i(\A)$ such that 
$\chi(\A^{X_{i+1}};t)$ divides $\chi(\A^{X_i} \setminus \{X_{i+1}\};t)$ for  
$i=0,\ldots,\ell-3$ and that 
$\A \setminus \{X_1\}$ is free.
\label{tower3}
\end{theorem}

\noindent
\textbf{Proof}. Follows immediately 
from Theorems \ref{newadd}, \ref{tower}, \ref{adddel} and 
Definitions \ref{DF0} and \ref{IF}. \owari
\medskip

\begin{cor}
$\A \in \mathcal{DF}$ if and only if 
$\A$ has a flag
$
V=X_0 \supset X_1 \supset \cdots \supset X_{\ell-2}$ with $X_i \in 
L_i(\A)$ such that 
$\A^{X_i}$ is all free for $i=1,\ldots,\ell-2$,
$\exp(\A^{X_{i+1}}) \subset \exp(\A^{X_i}) $ for  
$i=1,\ldots,\ell-3$ and that 
$\A \setminus \{X_1\}$ is free with $\exp(\A \setminus \{X_1\} )
\supset \exp(\A^{X_{1}})$.
\label{tower5}
\end{cor}

\noindent
\textbf{Proof}. 
Apply the same proof as that of Theorem \ref{tower3}.\owari
\medskip

The points in Theorem \ref{tower3} and Corollary \ref{tower5} are, though we have to check the 
freeness of some arrangements, we do not have to check 
anything else on them contrary to 
inductively free arrangements. 
Corollary \ref{tower5} is useful when we have a list of free arrangements.

In the end of this section let us prove Theorem \ref{b2free1}.
\medskip

\noindent
\textbf{Proof of Theorem \ref{b2free1}}. 
Let $\A$ be an arrangement and $H \in \A$. 
Recall that, if both $\A$ and 
$\A^H$ are free with $\exp(\A^H) \subset \exp(\A)$, then 
$\chi_0(\A;t)=(t-(|\A|-|\A^H|))\chi_0(\A^H;t)$ by 
the deletion theorem (Theorem \ref{adddel}) and the factorization theorem
(Theorem \ref{fact}). 
Hence  the implication (1) to (2) is immediate by the definition of divisional free. 
Assume (2). 
Let us agree that $X_{\ell-1}$ is an arrangemnet of one point on a line, we have 
\begin{eqnarray*}
b_2(d\A)&=&\sum_{0 \le i < j \le \ell-2}
(|\A^{X_i}|-|\A^{X_{i+1}}|)(|\A^{X_j}|-|\A^{X_{j+1}}|)\\
&=&
\sum_{i=0}^{\ell-3}
(|\A^{X_i}|-|\A^{X_{i+1}}|)\sum_{j=i+1}^{\ell-2}(|\A^{X_j}|-|\A^{X_{j+1}}|)\\
&=&
\sum_{i=0}^{\ell-3}
(|\A^{X_i}|-|\A^{X_{i+1}}|)(|\A^{X_{i+1}}|-1),
\end{eqnarray*}
which implies (3).

Finally, assume (3). 
Applying Theorem \ref{r0division} to $\A^{X_i}$ and $\A^{X_{i+1}}$, it holds that 
\begin{equation}
b_2(d\A^{X_i})\ge 
b_2(d\A^{X_{i+1}})+(|\A^{X_i}|-|\A^{X_{i+1}}|)
(|\A^{X_{i+1}}|-1).
\label{eq17}
\end{equation}
Taking the sum from $i=0$ to $i=\ell-3$ of the equation 
(\ref{eq17}), we have 
$$
b_2(d\A) \ge 
\sum_{i=0}^{\ell-3}(|\A^{X_i}|-|\A^{X_{i+1}}|)
(|\A^{X_{i+1}}|-1).
$$
Here we used the fact that $b_2(d\A^{X_{\ell-2}})=0$ since 
$\chi_0(\A^{X_{\ell-2}};t)=t-|\A^{X_{\ell-2}}|+1$.
Since (3) assumes the above inequality is the equality, all the inequalities in 
(\ref{eq17}) have to be the equalities. Now apply Theorem \ref{tower} to 
each $\A^{X_i}$ to show that they are all divisionally free, which completes the proof. \owari
\medskip

The following is an immediate corollary of the proof of Theorem \ref{b2free1}.

\begin{cor}
In the notation of Theorem \ref{b2free1}, let 
$V=X_0 \supset X_1 \supset \cdots \supset X_{\ell-2}$ be an arbitrary
flag with $X_i\in L_i(\A)$. Then 
$$b_2(d\A) \ge 
\sum_{i=0}^{\ell-3}(|\A^{X_i}|-|\A^{X_{i+1}}|)
(|\A^{X_{i+1}}|-1).$$
\label{flagineq}
\end{cor}

\section{Recursively and divisionally free arrangements}
First let us recall the definition of the recursively free arrangements.

\begin{define}[Recursively free arrangements]
A set of $\ell$-arrangements $\mathcal{RF}_\ell$ is defined by,
\begin{itemize}
\item[(a)] 
all arrangements are in $\mathcal{RF}_\ell$ when $\ell=1,2$, and 
\item[(b)]
for $\ell \ge 3$, $\A \in \mathcal{RF}_\ell$ if $\A=\emptyset_\ell$, or  
there is a hyperplane $H \in \A$ such that $\A^H\in \mathcal{RF}_{\ell-1},\ 
\A \setminus \{H\} \in \mathcal{RF}_\ell$ and 
$\chi(\A^H;t) \mid \chi(\A;t)$, or 
there is a hyperplane $L \not \in \A$ such that $(\A \cup \{L\}) \in 
\mathcal{RF}_{\ell},\ 
(\A \cup\{L\})^L \in 
\mathcal{RF}_{\ell-1}$ and $\chi((\A \cup\{L\})^L;t) \mid \chi(\A;t)$.
\end{itemize}
By Theorem \ref{adddel}, $\A$ is free if $\A \in \mathcal{RF}$. 
We say that an arrangement $\A$ is \textbf{recursively free} if 
$\A \in \mathcal{RF}:=\cup_{\ell \ge 1} \mathcal{RF}_\ell$. 
\label{RF0}
\end{define}

\begin{rem}
The definition of recursively free arrangements may seem to be artificial. 
However, this contains almost all of known free arrangements. There are only three arrangements known 
which are free but non-recursively free. The first one is found by Cuntz and Hoge in \cite{CH}, and 
the other two are in \cite{ACKN}. Though some free arrangements (like 
$\A(G_{31})$ in the previous section) are not known 
whether it is recursively free or not (see \cite{AHR}, Corollary 4.3), by these 
known results on $\mathcal{RF}$, 
it is worth considering Conjecture \ref{TC} for recursively free arrangements.
\label{RFTC}
\end{rem}

Contrary to the inductively free arrangements, it is not yet known whether 
the freeness of recursively free arrangements depends only on the combinatorics. 

\begin{problem}
Does the freeness of recursively free arrangements depend only on its 
combinatorics?
\label{pb}
\end{problem}

We will show that
a lot of recursively free arrangements are in fact divisionally free.
For that purpose, let us introduce a definition.

\begin{define}
Let $\mathcal{RF}=\cup_{\ell \ge 1}\mathcal{RF}_\ell$ denote the class of recursively free arrangements.
Define $\mathcal{RF}_1^*:=\mathcal{RF}_1,\ 
\mathcal{RF}_2^*:=\mathcal{RF}_2$, and 
\begin{eqnarray*}
\mathcal{RF}_\ell^*:&=&\{ \A \in \mathcal{RF}_\ell 
\mid \exists H  \in \A\ \mbox{such that}\ 
\chi(\A^H;t) \mid \chi(\A;t)\ \\
&\ &\mbox{and}\ 
\A^H \in \mathcal{RF}_{\ell-1}^*\}\ \ \ \ (\ell \ge 3),
\end{eqnarray*}
and put $\mathcal{RF}^*:=\cup_{\ell \ge 1} \mathcal{RF}^*_\ell$. We call 
$\mathcal{RF}^*$ the class of \textbf{recursively free arrangements at the 
hill}. If $\A \in \mathcal{RF} \setminus \mathcal{RF}^*$, then we say that 
$\A$ is a \textbf{recursively free arrangement at the valley}. 
\label{RFstar}
\end{define}

As far as we know, there have been no answer to Problem \ref{pb}. 
Here we give the first partial answer to Problem \ref{pb} asserting that 
recursively free arrangements not 
at the valley are divisionally free.

\begin{theorem}
$\mathcal{RF}^* \subset \mathcal{DF}$.
\label{nb}
\end{theorem}

\noindent
\textbf{Proof}.
We show by induction on $\ell$. When $\ell=1,2$, then there is nothing to show.
Assume that $\ell \ge 3$ and $\A \in \mathcal{RF}^*_\ell$. Then there is $ H \in \A$ such that 
$\chi(\A^H;t) \mid \chi(\A;t)$ and $\A^H \in \mathcal{RF}_{\ell-1}^*
\subset \mathcal{DF}_{\ell-1}$ by the induction hypothesis. 
Hence Definition \ref{DF} shows that $\A \in 
\mathcal{DF}_\ell$. \owari
\medskip

An example of non-inductively free $\A \in \mathcal{RF}^*$ is 
intermediate 
arrangements. Let 
$\A_\ell^k(r)\ (\ell \ge 2,\ 0 \le k \le \ell)$
be an $\ell$-arrangement defined by 
$$
\prod_{i=1}^k x_i \prod_{1 \le i<j\le \ell,\ 0 \le n <r}
(x_i-\zeta^n x_j )=0,
$$
where $\zeta$ is a primitive $r$-th root of unity. They are called \textbf{intermediate 
arrangements}, and first studied by 
Orlik and Solomon in \cite{OS}. They interpolate 
between the reflection arrangements 
corresponding to monomial groups, and they are all free arrangements (Propositions 
2.11 and 2.13, \cite{OS}). 
The inductive freeness of intermediate arrangements 
is studied in \cite{AHR}, and showed that 
$\A_\ell^k(r)$ is not inductively free if and only if 
$r \ge 3$ and $0 \le k \le \ell-3$ (Theorem 3.6, \cite{AHR}). On their 
divisional freeness, we have the following.

\begin{theorem}
Let $\ell \ge 3$. Then 
$\A_\ell^k(r) \in \mathcal{DF}$ if and only if $k \neq 0$. 
\label{akr}
\end{theorem}

\noindent
\textbf{Proof}. 
Let $H \in \A_\ell^k(r)$ be an arbitrary coordinate hyperplane. When 
$k =0$, by Corollary 1.3 in \cite{HR}, 
there are no $H \in \A$ such that $\A^H$ is free and $\chi(\A^H;t) \mid \chi(\A;t)$. 
Hence $\A_\ell^0(r)$ is not divisionally free. Assume that $k \neq 0$.  
Then for a coordinate hyperplane $H \in \A$, 
by the result in \cite{OS} (see Proposition 3.1 in \cite{AHR} too), 
it holds that 
$(\A_\ell^k(r))^H \simeq \A_{\ell-1}^{\ell-1}(r)$, and 
$$
\exp(\A_\ell^k(r)) \supset \exp(\A_\ell^k(r)^H),\ 
\exp(\A_j^j(r)) \supset \exp(\A_{j-1}^{j-1}(r)) \ \ (j=3,\ldots,\ell).
$$
Hence Theorem \ref{tower} shows that 
$\A_\ell^k(r) \in \mathcal{DF}$.\owari
\medskip

As shown in Theorem 3.6 in \cite{AHR}, 
$\A_\ell^k(r)\ (0 <k <\ell-2,\ \ell \ge 3)$ are not inductively free but 
recursively free. By Theorem \ref{akr}, 
the conditions in Theorem \ref{nb} are satisfied. Hence 
$\A_\ell^k(r) \in \mathcal{RF}^*$. So only the arrangement $\A_\ell^0(r)$ which is at the 
the valley of free paths of arrangements $\A_\ell^k(r)\ ( k\ge 0)$ is not 
divisionally free. 

Also, following the definition based on the addition-deletion theorems, 
we can define the following class of free arrangements.

\begin{define}
(1)\,\,
A set of $\ell$-arrangements $\mathcal{MF}_\ell$ is defined by,
\begin{itemize}
\item[(a)] 
all arrangements are in $\mathcal{MF}_\ell$ when $\ell=1,2$, and 
\item[(b)]
for $\ell \ge 3$, $\A \in \mathcal{MF}_\ell$ if 
there is a hyperplane $H \in \A$ such that $\A^H\in \mathcal{MF}_{\ell-1}$ and 
$\chi(\A^H;t) \mid \chi(\A;t)$, or 
there is a hyperplane $L \not \in \A$ such that $(\A \cup \{L\})^L \in 
\mathcal{MF}_{\ell-1}$ and $\chi((\A \cup\{L\})^L;t) \mid \chi(\A;t)$.
\end{itemize}
We say that an arrangement $\A$ is \textbf{multiplicatively free} if 
$\A \in \mathcal{MF}:=\cup_{\ell \ge 1} \mathcal{MF}_\ell$.

(2)\,\,
An $\ell$-arrangement is \textbf{hereditarily divisionally free} if 
$\A^X$ is divisionally free for all $X \in L(\A)$.
\label{MH}
\end{define}

Recall that the definition of hereditarily inductively free, where 
the definition of hereditarily divisionally free generalizes it.
An arrangement $\A$ is \textbf{hereditarily inductively free} if 
$\A^X$ is inductively free for all $X \in L(\A)$. Clearly, 
hereditarily inductively free arrangements are hereditarily divisionally free.

\begin{prop}
There is an arrangement $\A$ which is not hereditarily 
inductively free but hereditarily divisionally free. 
\label{notHIF}
\end{prop}

\noindent
\textbf{Proof}. 
Let us again use a divisionally free, but not inductively free unitary reflection arrangement 
$\A=\A(G_{31})$ in the proof of Theorem \ref{IFDF}. 
By Theorem 1.2 in \cite{AHR}, every $\A^X$ is inductively free 
for $X \in L_3(\A)$. Since $\A$ is a $4$-arrangement and 
$\mathcal{IF} \subsetneq \mathcal{DF}$ by Theorem \ref{DFTC}, the proof is 
completed.\owari

\begin{rem}
Apparently multiplicatively free arrangements are generalizations 
of recursively free arrangements. The former contains the latter. 
By the result in \cite{CH}, there is a 
free arrangement which is not multiplicatively free. Hence 
there is a free arrangement which is not either divisionally nor
multiplicatively free.
\end{rem}

By Theorem \ref{nb}, if there exists a counterexample to Conjecture 
\ref{TC} in $\mathcal{RF}$, then that has to be at the valley. So we may 
pose a problem.

\begin{problem}
Does the freeness of a recursively free arrangement at the valley 
depend only on its combinatorics?
\label{valley}
\end{problem}

Even if Problem \ref{valley} is settled positively, Conjecture \ref{TC} seems 
to be still difficult. However, Problem \ref{valley} is a natural one based on the results in this 
article. Based on several examples, we pose the following conjecture:

\begin{conj}
\begin{itemize}
\item[(1)]
$\mathcal{RF}^* \subsetneq \mathcal{DF}$.
\item[(2)]
$\mathcal{RF} \subsetneq \mathcal{MF}$.
\end{itemize}
\label{pbs}
\end{conj}

To settle Conjecture \ref{pbs} (1) positively, for example, 
recalling the proof of Theorem \ref{IFDF}, 
it 
suffices to show that $\A(G_{31})$ is not recursively free, which is 
not yet known (see \cite{AHR}, Corollary 4.3). It seems natural to believe that 
Conjecture \ref{pbs} (1) is true. On Conjecture \ref{pbs} (2), it seems still difficult 
to check since there are few known examples of a free arrangement which is non-recursively free 
(see \cite{CH}, \cite{ACKN}).

\section{Application to Weyl arrangements}

In this section assume that $\K=\R$ and consider a 
real irreducible crystallographic reflection arrangements and 
its deformations. 

First let us recall a notation on root systems and related Weyl arrangements.
Let $\Phi$ be a real irreducible crystallographic root system of rank $\ell$ with 
the Coxeter number $h$. 
Fix a positive system $\Phi^+$ of $\Phi$ and let 
$\Delta=\{\alpha_1,\ldots,\alpha_\ell\} \subset \Phi^+$ the associated simple system.
Let $W$ be the corresponding Weyl group and 
$H_\alpha$ the reflecting hyperplane with respect to the root $\alpha \in \Phi^+$ in 
$V=\R^\ell$. Then the \textbf{Weyl arrangement} 
$\A_W$ of $\Phi^+$ is defined by
$$
\A_W:=\{H_\alpha \mid \alpha \in \Phi^+\}.
$$
Also, the \textbf{$k$-extended Shi arrangement} $\Shi^k$ is defined by the equation 
$$
\Shi^k:=\{z=0\} \bigcup \{ H_\alpha^j\}_{\alpha \in \Phi^+,\ -k+1 \le j \le k}, 
$$
where $z$ is a new coordinate added to the vector space 
spanned by 
$\alpha_1,\ldots,\alpha_\ell$, and 
$H_\alpha^j:=\{\alpha=jz\}$ for $\alpha \in \Phi^+$ and $j \in \Z$.
The Shi arrangement was introduced by J.-Y Shi in \cite{Shi}, and 
has been well-studied by several mathematicians. See \cite{AT2} and 
\cite{Y1} for example.  
It is proved by Yoshinaga in \cite{Y1} that 
$\Shi^k$ is free with exponents $(1,kh,\ldots,kh)$.
We can show that 
$\Shi^k$ is divisionally free, which implies Theorem \ref{Shi} by Theorem \ref{DFTC}.

\begin{theorem}
$\Shi^k \in \mathcal{DF}$. 
\label{Shi2}
\end{theorem}

To prove Theorem \ref{Shi2}, we need a theorem.

\begin{theorem}
Let $\K$ be an arbitrary field, $V=\K^\ell$ and 
$\A$ be an
$\ell$-arrangement in $V$. 
Assume that $\A$ is free with $\exp(\A)=(1,d_1\ldots,d_{\ell-1})$, and that 
there are distinct hyperplanes 
$H_1,\ldots,H_{\ell-1} \in \A$ such that,
\begin{itemize}
\item[(1)]
$\A_i':=\A \setminus \{H_i\}$ is free with $\exp(\A_i')=
(1,d_1,d_2,\ldots,d_{i-1},d_i-1,d_{i+1},\ldots, d_{\ell-1})$ for $i=1,\ldots,\ell-1$.
\item[(2)]
$\A':=\A \setminus \{H_1,\ldots,H_{\ell-1}\}$ is free with 
$\exp(\A')=(1,d_1-1,d_2-1,\ldots,d_{\ell-1}-1)$.
\end{itemize}
%
Then $\A \in \mathcal{DF}$. 
\label{SRB}
\end{theorem}

\noindent
\textbf{Proof}. Let $X_i:=\cap_{j=\ell-i}^{\ell-1} H_j\ (i=1,\ldots,\ell-1)$. 
We show that 
$V=X_0 \supset X_1 \supset \cdots \supset X_{\ell-1}$ is a divisional flag.
Let $\alpha_j:=\alpha_{H_j}$ for $j=1,\ldots,\ell-1$.

For that purpose, at first, we show that $\A$ is free with $\exp(\A)=(1,d_1,d_2,\ldots,d_{\ell-1})$ with 
a special basis. By the assumption, there is a derivation $\theta_i \in D(\A_i')\ (i=1,\ldots,\ell-1)$ 
of degree $d_i-1$ such that $\theta_i \not \in D(\A)$, 
$\alpha_i \theta_i \in D(\A)$ and that $\theta_E,\theta_j$ are 
a part of a basis for $D(\A_j')$ for $j=1,\ldots,\ell-1$. We show that 
$\theta_E,\theta_1,\ldots,\theta_{\ell-1}$ form a basis for $D(\A')$. We may assume that 
$1 \le \deg \theta_1 \le \cdots \le \deg \theta_{\ell-1}$. 
Note that 
$1+\sum_{i=1}^{\ell-1} \deg \theta_i=|\A'|$ and $\deg \theta_j=d_j-1$ 
for $j=1,\ldots,\ell-1$. Hence by 
Proposition 4.42 in \cite{OT}, it suffices to show that 
$\langle \theta_E,\theta_1,\ldots,\theta_{i-1}\rangle_S \not \ni \theta_{i}$ for 
$i=1,\ldots,\ell-1$. If $\theta_1 \in S \theta_E$, then it contradicts the fact that 
$\theta_E,\theta_1$ form a part of basis for $D(\A_1')$. Now assume that 
$\theta_i=f\theta_E+\sum_{j=1}^{i-1} f_j \theta_j\ (f,f_j \in S)$. By the choice of $\theta_i$, 
$\theta_k(\alpha_j) \in S\alpha_j$ if and only if $k \neq j$. Since $\theta_E(\alpha)=\alpha$ for 
any $\alpha \in V^*$, it holds that $(f\theta_E+\sum_{j=1}^{i-1} f_j \theta_j)(\alpha_i) \in S\alpha_i$, but
$\theta_i(\alpha_i) \not \in S\alpha_i$, which is a contradiction. 

Hence $D(\A')$ has a basis $\theta_E,\theta_1,\ldots,\theta_{\ell-1}$ with the properties above. 
Moreover, the construction implies that $\theta_E,\alpha_1\theta_1,\ldots,\alpha_{\ell-1}\theta_{\ell-1}$ 
form a basis for $D(\A)$. By using these bases, 
we show that $V=X_0 \supset X_1 \supset \cdots \supset X_{\ell-1}$ is a divisional flag.

First, we show that $\A^{X_{1}}$ is free with $\exp(\A^{X_{1}})=(1,d_1,d_2,\ldots,d_{\ell-1})$. 
Let $\pi:D(\A) \rightarrow D(\A^{X_{1}})$ be the restriction map defined by 
$\pi(\theta)(f):=\theta(f)$ modulo $\alpha_{\ell-1}$. 
By Terao's restriction theorem \ref{rest}, the derivations 
$$\pi(\theta_E),\pi(\alpha_1\theta_1),\ldots,\pi(\alpha_{\ell-2}\theta_{\ell-2})$$ form a free basis for 
$D(\A^{X_{1}})$, hence $\exp(\A^{X_{1}})=(1,d_1,d_2,\ldots,d_{\ell-2})$. 

Second, the arrangement $\A^{X_{1}}$ in $H_{\ell-1}=\K^{\ell-1}$, 
hyperplanes $H_1 \cap H_{\ell-1},\ldots,H_{\ell-2} \cap H_{\ell-1} \in \A^{X_{1}}$ and 
derivations 
$\pi(\theta_E),\pi(\alpha_1\theta_1),\ldots,\pi(\alpha_{\ell-2}\theta_{\ell-2})$ satisfy the same assumptions as for 
$\A,H_1,H_2,\ldots,H_{\ell-1},\theta_E,\alpha_1 \theta_1,\ldots,\alpha_{\ell-1} \theta_{\ell-1}$. 
For $f \in S$, let 
$\overline{f}$ denote the image of $f$ by the canonical surjection $S \rightarrow S/(\alpha_{\ell-1})$. 
Since $\alpha_i \mid \theta_i$, clearly $\overline{\alpha_i} \mid \pi(\theta_i)$ 
for $i=1,\ldots,\ell-2$. Also, we can show that $\alpha_i=0\ (i=1,\ldots,\ell-2)$ 
define distinct hyperplanes in $\A^{X_{1}}$. Assume not. Then we may assume that 
$\overline{\alpha_2}=\overline{\alpha_3}$. Then for $\B_{1}:=
\A^{X_{1}} \setminus \{H_2 \cap H_{\ell-1}=H_3 \cap 
H_{\ell-1}\}$, the derivations 
$$
\pi(\theta_E),\pi(\alpha_1\theta_1),\pi(\theta_2),\pi(\theta_3), 
\pi(\alpha_4\theta_4),\ldots,\pi(\alpha_{\ell-2}\theta_{\ell-2})
$$
form a basis for $D(\B_{1})$. Since the sum of degrees of the derivations above 
is $|\A^{X_{1}}|-2=|\B_{1}|-1$, this contradicts Saito's criterion \ref{Saito}. Hence 
$H_1\cap H_{\ell-1},\cdots,H_{\ell-2} \cap H_{\ell-1}$ are distinct. Hence these satisfy the same assumption in the previous paragraph.

Now apply the same arguments inductively to $\A^{X_{1}},\A^{X_{2}},\cdots,\A^{X_{\ell-2}}$ with 
Theorem \ref{fact} to show that 
$\chi(\A^{X_{i+1}};t) \mid \chi(\A^{X_{i}};t)$ for $i=0,\ldots,\ell-3$. Hence Theorem \ref{tower} shows 
that $\A \in \mathcal{DF}$. \owari
\medskip

\noindent
\textbf{Proof of Theorem \ref{Shi2}}.
First recall the result on simple-root basis in \cite{AT}. It asserts that, there 
are derivations $\phi_1,\ldots,\phi_\ell \in D(\Shi^k)$ of degree $kh$ 
such that 
$\theta_E,\phi_1,\ldots,\phi_\ell$ form a homogeneous basis for $D(\Shi^k)$ and 
$(\alpha_i-kz) \mid \phi_i\ (i=1,\ldots,\ell)$. We call this basis a simple-root basis minus 
(see \cite{AT} for details). 

By the property of simple root basis, the Shi arrangement and 
hyperplanes $\alpha_i-kz=0\ (i=1,\ldots,\ell)$ clearly 
satisfy the assumptions in Theorem \ref{SRB}, which implies that 
$\Shi^k \in \mathcal{DF}$. \owari
\medskip

If we use some specific basis for the derivation module, then 
we can show the following corollary similar to Theorem \ref{SRB}.

\begin{cor}
Let $\K$ be an arbitrary field, $V=\K^\ell$ and 
$\A$ be a free 
$\ell$-arrangement in $V$ with $\exp(\A)=(1,d,d_2,\ldots,d_{\ell-1})$. 
Assume that  
there are distinct hyperplanes 
$H_2,\ldots,H_{\ell-1} \in \A$ and 
derivations 
$\theta_E, \varphi,\theta_2,\ldots,\theta_{\ell-1}$ for $D(\A)$ 
such that $\deg \varphi=d,\ 
\deg \theta_i=d_i\ (i=2,\ldots,\ell-1)$ and $\alpha_i:=\alpha_{H_i}$ divides $\theta_i$ for $i=2,\ldots,\ell-1$. 
Then $\A \in \mathcal{DF}$. 
\label{SRB2}
\end{cor}

\noindent
\textbf{Proof}. The proof is the same as that of Theorem \ref{SRB}. 
\owari
\medskip

Let 
$\Delta_j$ be the set of all positive roots of 
height $j$. 
By applying the same proofs as in Theorem \ref{SRB} with Corollary \ref{SRB2}, the following statements 
are easy to show by using results in \cite{ABCHT} and \cite{AT2}.

\begin{cor}
The following arrangements are divisionally free.
\begin{itemize}
\item[(1)]
$\A=\{H_\alpha^0\}_{\alpha \in \Delta_1} \cup 
\{H_\beta^0\}_{\beta \in \Sigma}$ for any subset 
$\Sigma \subset \Delta_2$.
\item[(2)]
$\A=\{H_\alpha^0\}_{\alpha \in \Delta_1 \cup \Delta_2 \cup \Delta_3}$.
\end{itemize}
\end{cor}

\noindent
\textbf{Proof}. 
(1)\,\,
By \cite{ABCHT}, 
$\A$ is free with basis $\theta_1,\ldots,\theta_{\ell-\sigma},\varphi_1,\ldots,\varphi_\sigma$ 
for $D(\A)$ 
with $\sigma:=|\Sigma|$ such that $\deg \theta_i=1\ (i=1,\ldots,\ell-\sigma),\ 
\deg \varphi_j=2\ (j=1,\ldots,\sigma)$ and that $\beta_j \mid \varphi_j$ for 
$j=1,\ldots,\sigma$ and $\Sigma_2=\{\beta_1,\ldots,\beta_\sigma\}$. By the argument in 
the proof of Theorem \ref{SRB}, we can construct a divisional flag
$V=X_0 \supset X_1 \supset \cdots \supset X_{\sigma}=X$, where 
$X_i=\cap_{j=1}^i H_{\beta_j}$. 
Since $\A^X$ is a Boolean arrangement, 
it is divisionally free. So $\A$ is divisionally free.

(2)\,\, By \cite{ABCHT}, $\A$ has a free basis for $D(\A)$ 
satisfying the condition in Corollary \ref{SRB2}. Hence $\A$ is divisionally free.\owari
\medskip

\begin{cor}
The following arrangements are divisionally free.
\begin{itemize}
\item[(1)]
$\Shi^k \setminus \{H_\alpha^k\}$ for $\alpha \in \Delta_1$.
\item[(2)]
$\Shi^k \setminus \{H_\alpha^k\}_{\alpha \in \Delta_1}$.
\item[(3)]
$\Shi^k \cup \{H_\alpha^{-k} \}_{\alpha \in \Delta_1}$.
\item[(4)]
$\Shi^k \cup \{H_\alpha^{-k} \}_{\alpha \in \Delta_1} \setminus 
\{H_\beta^{-k}\}$ for $\beta \in \Delta_1$.
\item[(5)]
$\Shi^k \cup \{H_\alpha^{-k} \}_{\alpha \in \Delta_1 \cup \Delta_2} $.
\end{itemize}
\end{cor}

\noindent
\textbf{Proof}. It was shown in \cite{AT2} that all the above cases have 
the free bases satisfying the condition in Corollary \ref{SRB2}. Hence they are 
all divisionally free. 
\owari
\medskip



Based on the arguments in this section, 
\cite{ABCHT}, \cite{AT2} and \cite{Y1}, we pose 
a conjecture. For that purpose, let us recall terminologies. For 
$\alpha,\beta \in \Phi^+$ For the terminology, 
we say $\alpha \ge \beta$ if $\alpha -\beta$ is a linear combination of 
simple roots with non-negative coefficients. An \textbf{ideal} $I \subset \Phi^+$ 
is defined by, if $\alpha, \beta \in \Phi^+$ with $\alpha \ge \beta$ and 
$\alpha \in I$ implies $\beta \in I$. Then an \textbf{ideal-subarrangement} $\A_I \subset \A_W$ 
is defined by 
$\A_I:=\{H_\alpha^0 \mid \alpha \in I\}$. Also, \textbf{ideal-Shi arrangement} 
$\Shi^k_{\pm I}$ is defined by 
\begin{eqnarray*}
\Shi^k_{+I}:&=&\Shi^k \cup \{H_\alpha^{-k}\}_{\alpha \in I},\\
\Shi^k_{-I}:&=&\Shi^k \setminus \{H_\alpha^{k}\}_{\alpha \in I}.
\end{eqnarray*}
For details, 
see \cite{ABCHT}
and \cite{AT2}. 

\begin{conj}
All ideal-subarrangments $\A_I$ and ideal-Shi arrangements $\Shi^k_{\pm I}$ are 
divisionally free.
\label{ideals}
\end{conj}

\noindent
\textbf{Acknowledgements}. The author is partially supported by 
by JSPS Grants-in-Aid for Young Scientists
(B)
No. 24740012.
Also, the author is grateful to Torsten Hoge for his pointing out an error 
in Theorem 6.2.

\end{document}